\documentclass[11pt]{article}
\title{On incidence algebras description of cobweb
posets}
\author{Ewa Krot-Sieniawska \\
\\Institute of Computer Science, Bia{\l}ystok University\\
PL-15-887 Bia{\l}ystok, ul.Sosnowa 64, POLAND\\
e-mail: ewakrot@wp.pl, ewakrot@ii.uwb.edu.pl}
\date{}

\usepackage{amsmath,amsthm}
\usepackage{amssymb}
\usepackage{graphicx}
\chardef\bslash=`\\ 
\newtheorem{defn}{Definition}[]

\newtheorem{rem}{Remark}[]
\newtheorem{thm}{Theorem}[]

\begin{document}
\maketitle
\begin{abstract}
\noindent The explicite formulas for M\"{o}bius function and some
other important elements of the incidence algebra of an arbitrary
cobweb poset are delivered.  For that to do one uses
Kwa\'sniewski's construction of his  cobweb posets \cite{44,46}.
The digraph representation of these cobweb posets
constitutes a newly discovered class of orderable DAG's \cite{p,eksKoDAG,MD}
named here down KoDAGs  with a kind of universality now being investigated.
Namely  cobweb posets'  and thus KoDAGs's defining
di-bicliques are links of any complete relations' chains.

\end{abstract}
AMS Classification numbers: 06A06,  06A07, 06A11,  11C08, 11B37\\

\noindent Presented at Gian-Carlo Rota Polish Seminar:
http://ii.uwb.edu.pl/akk/sem/sem rota.htm
\section{ Cobweb posets}
The   family of the so called cobweb posets $\Pi$ has been
invented by A.K.Kwa\'sniewski few years ago (for references  see:
\cite{44,46}). These structures are such a generalization of the
Fibonacci tree growth  that allows joint combinatorial
interpretation for all of them under the admissibility condition
(see \cite{49,49a}).

\noindent Let $\{F_n\}_{n\geq 0}$ be a natural numbers valued
sequence with $F_0=1$ (with $F_0=0$ being exceptional as in case
of Fibonacci numbers). Any sequence satisfying this property
uniquely designates cobweb poset  defined as follows.

\noindent For $s\in\bf{N}_0=\bf{N}\cup\{0\}$ let us to define
levels of
 $\Pi$:
$$\Phi_{s}=\left\{\langle j,s \rangle ,\;\;1\leq j \leq F_{s}\right\},\;\;\;$$
(in case of $F_0=0$  level $\Phi_0$ corresponds to the empty root
$\{\emptyset\}$). )

\noindent Then

\begin{defn}
Corresponding cobweb poset is  an infinite partially ordered set
$\Pi=(V,\leq)$, where
 $$ V=\bigcup_{0\leq s}\Phi_s$$
 are the elements ( vertices) of $\Pi$ and the partial order relation $\leq$ on $V$ for
 $x=\langle s,t\rangle, y=\langle u,v\rangle $ being  elements of
cobweb poset $\Pi$ is defined by  formula
$$ ( x \leq_{P} y) \Longleftrightarrow
 [(t<v)\vee (t=v \wedge s=u)].$$
\end{defn}

\noindent Obviously any cobweb poset can be represented, via its
Hasse diagram, as infinite directed  graf  $\Pi=\left(
V,E\right)$, where  set $V$ of its vertices is defined as above
and

$$E =\{\left(\langle j , p\rangle,\langle q ,(p+1) \rangle
\right)\}\;\cup\;\{\left(\langle 1 , 0\rangle ,\langle 1 ,1
\rangle \right)\},$$ \quad where $1 \leq j \leq {F_p}$ and $1\leq
q \leq {F_{(p+1)}}$ stays for  set of (directed) edges.





\noindent The Kwasniewski cobweb posets under consideration
represented by  graphs  are examples of oderable directed acyclic
graphs (oDAG)  which we start to call from now in brief:  KoDAGs.
These are  structures of universal importance for the whole of
mathematics - in particular for discrete "`mathemagics"'
[http://ii.uwb.edu.pl/akk/ ]  and computer sciences in general
(quotation from \cite{49,49a} ):

\begin{quote}
For any given natural numbers valued sequence the graded (layered)
cobweb posets` DAGs  are equivalently representations of a chain
of binary relations. Every relation of the cobweb poset chain is
biunivocally represented by the uniquely designated
\textbf{complete} bipartite digraph-a digraph which is a
di-biclique  designated  by the very  given sequence. The cobweb
poset is then to be identified with a chain of di-bicliques i.e.
by definition - a chain of complete bipartite one direction
digraphs.   Any chain of relations is therefore obtainable from
the cobweb poset chain of complete relations  via
deleting  arcs (arrows) in di-bicliques.\\
Let us underline it again : \textit{any chain of relations is
obtainable from the cobweb poset chain of  complete relations via
deleting arcs in di-bicliques of the complete relations chain.}
For that to see note that any relation  $R_k$ as a subset of  $A_k
\times A_{k+1}$ is represented by a  one-direction bipartite
digraph  $D_k$.  A "complete relation"  $C_k$ by definition is
identified with its one direction di-biclique graph $d-B_k$.  Any
$R_k$ is a subset of  $C_k$. Correspondingly one direction digraph
$D_k$ is a subgraph of an one direction digraph of $d-B_k$.\\
The one direction digraph of  $d-B_k$ is called since now on
\textbf{the di-biclique }i.e. by definition - a complete bipartite
one direction digraph.   Another words: cobweb poset defining
di-bicliques are links of a complete relations' chain.
\end{quote}

\noindent According to the  definition above arbitrary cobweb
poset $\Pi=(V,\leq)$ is a graded poset ( ranked poset) and for
$s\in\bf
 N_0$:
 $$x\in\Phi_s\;\; \longrightarrow\;\; r(x)=s,$$
 where $r:\Pi\rightarrow \bf N_0$ is a rank function on $\Pi$.

\noindent Let us then define Kwa\'sniewski finite cobweb
sub-posets as follows
 \begin{defn}
 Let $P_n=(V_n,\leq)$, $(n\geq 0)$, for ${\displaystyle
V_n=\bigcup_{0\leq s\leq n}}\Phi_s$ and  $\leq$ being the induced
partial order relation on $\Pi$.
\end{defn}

\noindent Its easy to see that $P_n$ is ranked poset with rank
function $r$ as above. $P_n$ has a unique minimal element
$0=\langle 1,0\rangle$ ( with $r(0)=0$). Moreover $\Pi$ and all
$P_n$ s are locally finite, i.e. for any pair $x,y\in \Pi$, the
segment $[x,y]=\{z\in\Pi:\,x\leq z\leq y\}$ is finite.

 I this paper we shall consider the incidence algebra of an
 arbitrary cobweb poset. The one for Fibonnaci cobweb poset
 uniquely designated by the famous Fibonacci sequence was
 presented by the present author in \cite{35, 36} where M\"{o}bius function and some other
important elements of the incidence algebra were delivered. As we
shall see, the construction of an arbitrary cobweb poset universal
for all such structures enables us to extend these results to the
hole family of cobweb posets.

\section{Incidence algebra of an arbitrary cobweb poset}

Let us recall that one defines the incidence algebra of a locally
finite partially ordered set $P$ as follows (see \cite{65, 73,
74}):
$$ { I(P)}=I(P,R)=\{f: P\times  P\longrightarrow  R;\;\;\;\; f(x,y)=0\;\;\;unless\;\;\; x\leq y\}.$$
The sum of two such functions $f$ and $g$ and multiplication by
scalar are defined as usual. The product $h=f\ast g$ is defined as
follows:
$$ h(x,y)=(f\ast g)(x,y)=\sum_{z\in {\bf P}:\;x\leq z\leq y} f(x,z)\cdot g(z,y).$$
It is immediately verified that this is an associative algebra
(with  an identity element $\delta (x,y)$, the Kronecker delta),
over  any associative ring  R.

Let $\Pi$ be an arbitrary cobweb poset uniquely designated
 by the natural numbers valued sequence $\{F_n\}_{n\geq 0}$, in the way  as written above. Now,  we shall construct
 some typical elements of incidence algebra   ${I(\Pi)}$ of $\Pi$.
 Let $x,y$ be some arbitrary elements of $\Pi$ such that $x=\langle s,t\rangle$, $y=\langle u,v\rangle
  $, $(s,u\in\bf{N}$, $t,v\in \bf{N}_0)$, $1\leq s\leq F_t$
  and $1\leq u\leq F_v$.

 The zeta function of ${\Pi}$  defined by:
$$ \zeta (x,y)=\Big\{\begin{array}{l}1\;\;for\;\;\; x \leq y\\0\;\;otherwise\end{array}$$
 is an element of ${I(\Pi)}$. Obviously it is a characteristic
 function of partial order in $\Pi$. One can show that for $x,y\in
 \Pi$ as above
 \begin{equation}
\zeta(x,y)=\zeta\left( \langle s,t\rangle ,\langle u,v \rangle
\right)=\delta (s,u)\delta (t,v)+\sum_{k=1}^{\infty}\delta(t+k,v),
\end{equation}

The knowledge of $\zeta$ enables us to construct other typical
elements of incidence algebra  of ${\Pi}$. The one of them is the
M\"{o}bius function indispensable in numerous inversion type
formulas of countless applications \cite{65, 73, 74}. Of course
the $\zeta$ function of a locally finite partially ordered set is
invertible in incidence algebra and its inversion is the famous
M\"{o}bius function $\mu$ i.e.:
$$\zeta \ast \mu=\mu \ast \zeta=\delta.$$
One can recover it just by the use of the recurrence formula for
M\"{o}bius function of locally finite partially ordered set
${I(P)}$ (see \cite{65}):
\begin{equation}
 \left\{\begin{array}{l}\mu(x,x)=1\;\;\;\;for\;\; all \;\;x\in{\bf
 P}\quad \quad \quad \\\ \\
\mu (x,y)=-\sum_{x\leq z<y}\mu (x,z)\end{array}\right.
\end{equation}
Namely for $x,y\in\Pi$ as above one has
\begin{multline}\label{mobius2}
\mu(x,y)=\mu\left( \langle s,t\rangle ,\langle u,v\rangle
\right)=\\
\qquad \quad \;
=\delta(t,v)\delta(s,u)-\delta(t+1,v)+\sum_{k=2}^{\infty}\delta(t+k,v)(-1)^{k}
\prod_{i=t+1}^{v-1}(F_{i}-1),
\end{multline}
for
$$ \delta (x,y)=\Big\{\begin{array}{l}1\;\;\;\;\; x=y\\0\;\;\;\;\;x\neq
y\end{array}.$$

The formula (\ref{mobius2}) enables us to formulate the following
theorem (see \cite{65}):
\begin{thm}{\bf (M\"{o}bius Inversion Formula for $\Pi$)}\\
Let $f(x)=f(\langle s,t\rangle)$ be a $R$ valued function, defined
for $x=\langle s,t\rangle$ ranging in cobweb poset $\Pi$. Let an
element $p=\langle p_{1},p_{2}\rangle$ exist with the property
that $f(x)=0$ unless $x\geq p$.

Suppose that $$g(x)=\sum_{\left\{ y \in P :\; y\leq x
\right\}}f(y).$$ Then
$$f(x)=\sum_{\{y\in P:\;y\leq x\}}g(y)\mu(y,x).$$
Hence using coordinates of $x,y$ in $\Pi$ i.e. $x=\langle
s,t\rangle,\;y=\langle u,v\rangle$ if
$$g(\langle s,t\rangle)=\sum_{v=0}^{t-1}\sum_{u=1}^{F_{v}}\left(f(\langle u,v\rangle )\right)+f(\langle s,t\rangle)$$
then we have
\begin{multline}
f(\langle s,t\rangle)=\sum_{v\geq 0}\sum_{u=1}^{F_{v}}g(\langle
u,v\rangle)\mu(\langle s,t\rangle ,\langle
u,v\rangle)=\\
=\sum_{v\geq 0}\sum_{u=1}^{F_{v}}g(\langle
u,v\rangle)\left[\delta(v,t)\delta(u,s)-\delta(v+1,t)+\sum_{k=2}^{\infty}\delta(v+k,t)(-1)^{k}
\prod_{i=v+1}^{t-1}(F_{i}-1)\right].
\end{multline}
\end{thm}

 Now we shall  deliver some  other typical
elements of incidence algebra $I(\Pi)$ perfectly suitable for
calculating number of chains, of maximal chains etc. in finite
sub-posets of ${\Pi}$.

The function $\zeta ^{2}=\zeta \ast \zeta$ counts the number of
elements in the segment $\left[ x,y\right]$ (where $x=\langle s,t
\rangle , y=\langle u,v \rangle$), i.e.:
\begin{multline*}
\zeta^{2}(x,y)=\left(\zeta \ast \zeta\right) (x,y)=\sum_{x \leq z
\leq y}\zeta (x,z)\cdot\zeta(z, y)=\sum_{x \leq z \leq
y}1=\mathrm{card} \left[ x, y\right]
\end{multline*}
Therefore for $x,y\in \Pi$ as above, we have: \begin{equation}
\begin{array}{lll}
  card \left[ x,
y\right] & = &
\left({\displaystyle\sum_{i=t+1}^{v-1}\sum_{j=1}^{F_{i}}}1\right)+2
=  \left({\displaystyle\sum_{i=t+1}^{v-1}}F_{i}\right)+2 .\\
\end{array}\end{equation}

For any incidence algebra the function $\eta$ is defined as
follows: $$\eta(x,y)=(\zeta-\delta)(x,y)=\left\{
\begin{array}{lll}
1&&x<y\\
0&&otherwise \end{array}\right.$$ The corresponding function for
$x,\; y$ being elements of the cobweb poset ${\Pi}$, ($x=\langle
s,t \rangle , y=\langle u,v \rangle$) is then given by formula:
\begin{equation}
\eta(x,y)=\sum_{k=1}^{\infty}\delta(t+k,v)=\left\{\begin{array}{ccc}
  1 &  & t<v \\
  0 &  & w\;p.p. \\
\end{array}\right..\end{equation}
It was shown \cite{65, 74} that $\eta ^{k}(x,y),\;\; (k\in {\bf
N})$ counts the number of chains of length $k$, (with $(k+1)$
elements) from $x$ to $y$. In $\Pi$ one has
\begin{equation}\begin{array}{lll}
   \eta^2(x,y) & = & {\displaystyle \sum_{x\leq z\leq y}}\eta(x,z)\eta(z,y) \\
    & = & {\displaystyle \sum_{x< z< y}}1 =\mathrm{card}[x,y]-2 = F_{v+1}-F_{t+2}, \\
 \end{array}\end{equation}
(for $F_{v+1}-F_{t+2}<0$  one takes $\eta^2(x,y)=0$) and
\begin{equation}\begin{array}{lll}
   \eta^3(x,y) & = & {\displaystyle \sum_{x\leq z_1\leq z_2\leq y}}\eta(x,z_1)\eta(z_1,z_2)\eta(z_2,y) \\
    & = & {\displaystyle \sum_{x< z_1<z_2< y}}1 ={\displaystyle \sum_{t<k<l<v}} F_{k}F_{l}.\\
 \end{array}\end{equation}
In general, for $k\geq 0$:
\begin{equation}
\begin{array}{lll}
   \eta^k(x,y) & = & {\displaystyle \sum_{x< z_1<z_2<...<z_{k-1}<
   y}}1\\&&\\
   & =&{\displaystyle \sum_{t<i_1<i_2<...<i_{k-1}<v}} F_{i_1}F_{i_2}...F_{i_{k-1}}.\\
 \end{array}\end{equation}

 Now let
$$\mathcal{C}(x,y)=(2\delta-\zeta)(x,y)=\left\{\begin{array}{lll}
1&&x=y\\
-1&&x<y\\
0&&otherwise \end{array} \right.$$ For elements of ${\Pi}$ as
above we have:
\begin{equation}\mathcal{C}(\langle s,t\rangle , \langle
u,v\rangle)=
\delta(t,v)\delta(s,u)-\sum_{k=1}^{\infty}\delta(t+k,v).\end{equation}
The inverse function $\mathcal{C}^{-1}(x,y)$ counts the number of
all chains from $x$ to $y$. From the recurrence formula one infers
that
$$\left\{\begin{array}{l}
\mathcal{C}^{-1}(x,x)=\frac{1}{\mathcal{C}(x,x)}\\
\\
\mathcal{C}^{-1}(x,y)=-\frac{1}{\mathcal{C}(x,x)}\sum_{x< z\leq
y}\mathcal{C}(x,z)\cdot\mathcal{C}^{-1}(z,y) \end{array}\right.$$

 For any incidence algebra the function $\chi$ is defined as
follows:
$$\chi (x,y)=\left\{\begin{array}{lll}
1&&x \lessdot y\\
0&& w\,p.p., \end{array}\right.,
 $$
 where $x\lessdot y$ iff $y$ covers $x$, i.e. $|[x,y]|=2$. For $x,y\in\Pi$ as above one has
\begin{equation} \chi(x,y)=\delta(t+1,v).\end{equation}
It was shown \cite{65,73,74} that $\chi ^{k}(x,y),\;\; (k\in {\bf
N})$ counts the number of maximal chains of length $k$, (with
$(k+1)$ elements) from $x$ to $y$.  In $\Pi$ one has
\begin{equation}
\chi^2(x,y)=\sum_{x\lessdot z\lessdot y}1=\delta (t+2,v)F_{t+1},
\end{equation}
and
\begin{equation}
\chi^3(x,y)=\sum_{x\lessdot z_1\lessdot z_2 \lessdot y}1=\delta
(t+3,v)F_{t+1}F_{t+2}.
\end{equation}
In general
\begin{equation}
\chi^k(x,y)=\sum_{x\lessdot z_1\lessdot ...\lessdot
z_{k-1}\lessdot y}1=\delta (t+k,v)F_{t+1}F_{t+2}...F_{v-1}.
\end{equation}
for $k\geq 0$.

 Finally let
 $$\mathcal{M}(x,y)=(\delta -\chi)(x,y)=\left\{
\begin{array}{ll}
1&x=y\\
-1&x \lessdot y\\
0&otherwise \end{array}\right.$$
 For elements of ${\Pi}$, ($x=\langle
s,t \rangle , y=\langle u,v \rangle$) one has:
\begin{equation}\mathcal{M}(\langle s,t\rangle ,\langle
u,v\rangle)=\delta(t,v)\delta(s,u)-\delta(t+1,v).\end{equation}
Then the  inverse function of $\mathcal{M}$:
$$\mathcal{M}^{-1}=\frac{\delta}{\delta -\chi}=\delta +\chi +\chi
^{2}+\chi^{3}+\ldots$$ counts the number of all maximal chains
from $x$ to $y$.

\begin{rem}{\em
Let us remark that above formulas hold for  functions: $\zeta$,
$\mu$, $\eta$, etc.   being  elements of $I(P_n)$ for $n\geq 0$ -
the incidence algebras of finite cobweb posets $P_n$ i.e. finite
sub-posets of $\Pi$. For example in arbitrary $P_n$ one has
 \begin{equation}
\zeta(x,y)=\zeta\left( \langle s,t\rangle ,\langle u,v \rangle
\right)=\delta (s,u)\delta (t,v)+\sum_{k=1}^{n}\delta(t+k,v),
\end{equation}
for $x=\langle s,t\rangle$, $y=\langle u,v\rangle
  $, $0\leq t,v\leq n$, $1\leq s\leq F_t$
  and $1\leq u\leq F_v$. }\end{rem}
  \begin{rem}{\em
 Let us recall that for $P$ being finite poset, the incidence
 algebra $I(P)$ over a commutative ring R with identity is
 isomorphic to a subring of $M_{|P|}(R)$, i.e. subring of all
 upper triangular $|P|\times |P|$ matrices over the ring $R$, \cite{65, 73,74}.

 Let us show how it works in case of (finite) cobweb poset $P_n$, ($n\geq 0$).
 One can define a chain (linear order) $X=(X,\leq_X)$ on the set of all elements of $P_n$
 as follows:
 $$ ( \langle s,t\rangle \leq_{{ X}} \langle u,v\rangle) \Longleftrightarrow
 [(t < v)\vee (t=v \wedge s\leq u)].$$
 Now for $f$ being an element of $I(P_n)$ let us define
 corresponding  matrix
 $M(f)=[m_{ij}]$ as follows
$$m_{i,\,j}=f(x_i,x_j),$$
where  $x_i,x_j$ are $i$-th i $j$-th elements  of the chain $X$,
respectively. It's easy to verify that that $M(f)$ is a
$\nu\times\nu$, ($\nu=1+\sum_{k=1}^nF_k$) upper triangular matrix.
Then the product $f\ast g$ corresponds to the product $M(f)M(g)$
of matrices and an invertible element $f\in I(p_n)$ corresponds to
an invertible matrix $M(f)$, i.e. $\det M(f)\neq 0$.

 For example matrix $M(\zeta)$ corresponding to $\zeta\in P_6$
 being a finite cobweb poset designated by the sequence of Fibonacci numbers
 is of the form
 \begin{center}
 {\footnotesize  $
M(\zeta)=\left[\begin{array}{ccccccccccccccccccccc}
  1 & 1 & 1 & 1 & 1 & 1 & 1 & 1 & 1 & 1 & 1 & 1 & 1 & 1 & 1 & 1 & 1 & 1 & 1 & 1 & 1 \\
  0 & 1 & 1 & 1 & 1 & 1 & 1 & 1 & 1 & 1 & 1 & 1 & 1 & 1 & 1 & 1 & 1 & 1 & 1 & 1 & 1 \\
  0 & 0 & 1 & 1 & 1 & 1 & 1 & 1 & 1 & 1 & 1 & 1 & 1 & 1 & 1 & 1 & 1 & 1 & 1 & 1 & 1 \\
  0 & 0 & 0 & 1 & 0 & 1 & 1 & 1 & 1 & 1 & 1 & 1 & 1 & 1 & 1 & 1 & 1 & 1 & 1 & 1 & 1 \\
  0 & 0 & 0 & 0 & 1 & 1 & 1 & 1 & 1 & 1 & 1 & 1 & 1 & 1 & 1 & 1 & 1 & 1 & 1 & 1 & 1 \\
  0 & 0 & 0 & 0 & 0 & 1 & 0 & 0 & 1 & 1 & 1 & 1 & 1 & 1 & 1 & 1 & 1 & 1 & 1 & 1 & 1 \\
  0 & 0 & 0 & 0 & 0 & 0 & 1 & 0 & 1 & 1 & 1 & 1 & 1 & 1 & 1 & 1 & 1 & 1 & 1 & 1 & 1 \\
  0 & 0 & 0 & 0 & 0 & 0 & 0 & 1 & 1 & 1 & 1 & 1 & 1 & 1 & 1 & 1 & 1 & 1 & 1 & 1 & 1 \\
  0 & 0 & 0 & 0 & 0 & 0 & 0 & 0 & 1 & 0 & 0 & 0 & 0 & 1 & 1 & 1 & 1 & 1 & 1 & 1 & 1 \\
  0 & 0 & 0 & 0 & 0 & 0 & 0 & 0 & 0 & 1 & 0 & 0 & 0 & 1 & 1 & 1 & 1 & 1 & 1 & 1 & 1 \\
  0 & 0 & 0 & 0 & 0 & 0 & 0 & 0 & 0 & 0 & 1 & 0 & 0 & 1 & 1 & 1 & 1 & 1 & 1 & 1 & 1 \\
  0 & 0 & 0 & 0 & 0 & 0 & 0 & 0 & 0 & 0 & 0 & 1 & 0 & 1 & 1 & 1 & 1 & 1 & 1 & 1 & 1 \\
  0 & 0 & 0 & 0 & 0 & 0 & 0 & 0 & 0 & 0 & 0 & 0 & 1 & 1 & 1 & 1 & 1 & 1 & 1 & 1 & 1 \\
  0 & 0 & 0 & 0 & 0 & 0 & 0 & 0 & 0 & 0 & 0 & 0 & 0 & 1 & 0 & 0 & 0 & 0 & 0 & 0 & 0 \\
  0 & 0 & 0 & 0 & 0 & 0 & 0 & 0 & 0 & 0 & 0 & 0 & 0 & 0 & 1 & 0 & 0 & 0 & 0 & 0 & 0 \\
  0 & 0 & 0 & 0 & 0 & 0 & 0 & 0 & 0 & 0 & 0 & 0 & 0 & 0 & 0 & 1 & 0 & 0 & 0 & 0 & 0 \\
  0 & 0 & 0 & 0 & 0 & 0 & 0 & 0 & 0 & 0 & 0 & 0 & 0 & 0 & 0 & 0 & 1 & 0 & 0 & 0 & 0 \\
  0 & 0 & 0 & 0 & 0 & 0 & 0 & 0 & 0 & 0 & 0 & 0 & 0 & 0 & 0 & 0 & 0 & 1 & 0 & 0 & 0 \\
  0 & 0 & 0 & 0 & 0 & 0 & 0 & 0 & 0 & 0 & 0 & 0 & 0 & 0 & 0 & 0 & 0 & 0 & 1 & 0 & 0 \\
  0 & 0 & 0 & 0 & 0 & 0 & 0 & 0 & 0 & 0 & 0 & 0 & 0 & 0 & 0 & 0 & 0 & 0 & 0 & 1 & 0 \\
  0 & 0 & 0 & 0 & 0 & 0 & 0 & 0 & 0 & 0 & 0 & 0 & 0 & 0 & 0 & 0 & 0 & 0 & 0 & 0 & 1 \\
\end{array}\right]$}
\end{center}
 The formulas delivered above allow us to construct matrices of
 this and  other elements of $I(P_n)$ for $P_n$ being the finite
 cobweb poset designated by an arbitrary natural number valued
 sequence $\{F_n\}$ with $F_0=1$ ($F_0=0$ being exceptionable as in case of Fibonacci
 sequence).

  }\end{rem}

\noindent {\bf Acknowledgements}\\
Discussions with Participants of Gian-Carlo Rota Polish Seminar,\\
http://ii.uwb.edu.pl/akk/sem/sem\_rota.htm are highly appreciated.

 \end{document}